\documentclass[11pt]{article}

\usepackage{amstext}
\usepackage{amssymb}
\usepackage{amsmath}
\usepackage{pb-diagram} \usepackage{amsthm}

 \usepackage [latin1]{inputenc}
 
\usepackage{rotating}
\usepackage{tikz-cd}

\usepackage{amsfonts}
\usepackage{graphicx}
\usepackage[skip=0pt]{caption}
\theoremstyle{plain}
 \newtheorem{thm}{Theorem}
 \newtheorem{lem}{Lemma}

\newtheorem{defn}{Definition}

\newcommand{\noi}{\noindent}

\newenvironment{dedication}

\title{The number of Simply-connected Trivalent $2$-dimensional Stratifolds}

\author{J. C. G\'{o}mez-Larra\~{n}aga\thanks{Centro de
Investigaci\'{o}n en Matem\'{a}ticas, A.P. 402, Guanajuato 36000, Gto. M\'{e}xico. jcarlos@cimat.mx} \and F.
Gonz\'alez-Acu\~na\thanks{Instituto de Matem\'aticas, UNAM, 62210 Cuernavaca, Morelos,
M\'{e}xico and Centro de
Investigaci\'{o}n en Matem\'{a}ticas, A.P. 402, Guanajuato 36000,
Gto. M\'{e}xico. fico@math.unam.mx} \and Wolfgang
Heil\thanks{Department of Mathematics, Florida State University,
Tallahasee, FL 32306, USA. heil@math.fsu.edu}}
\date{}

\begin{document}

\maketitle

\begin{dedication}\begin{center} Dedicated to Sergey Antonyan on the occasion of his 65th birthday
\end{center}       \vspace{\baselineskip}
      \end{dedication}

\begin{abstract} We describe a method for counting the number of $1$-connected trivalent $2$-stratifolds with a given number of singular curves and $2$-manifold components.
\end{abstract}



\section{Introduction}   

Observations in data analysis suggest that the points in a naturally-occurring dataset tend to cluster near a manifold with singularities. In particular, for dimension $2$, these manifolds with singularities are $2$-stratifolds and occur in the study of the energy landscape of cyclo-octane \cite{Martin}, with a systematic application of local topological methods described in \cite{Stolz}, the study of boundary singularities produced by the motion of soap films \cite{Gold}, and in organizing  data \cite{B2}. A systematic study of trivalent $2$-stratifolds was begun in \cite{GGH}. Whereas closed $2$-manifolds are classified by their fundamental groups, this is far from true for $2$-stratifolds. In fact, for any given $2$-stratifold there are infinitely many others with the same fundamental group. The question arises whether one can effectively construct all of the $2$-stratifolds that have a given fundamental group.

A $2$-stratifold is essentially determined by its associated bi-colored labeled graph  and a presentation for its fundamental group can be read off from the labeled graph. Thus the question arises when a labeled graph determines a simply connected $2$-stratifold. In \cite{GGH} an algorithm on the labeled graph was developed for determining whether the graph determines a simply connected $2$-stratifold and in \cite{GGH1} we obtained a complete classification of all trivalent labeled graphs that represent simply connected $2$-stratifolds. Then in \cite{GGH2}
we developed three operations on labeled graphs that will construct recursively from a single vertex all trivalent graphs that represent $1$-connected $2$-stratifolds. A referee of that paper asked whether it is possible to compute the number of all such labeled graphs for a given number of vertices. The purpose of the present paper is to describe a method that leads to such computations. Our approach is based on the classification theorem in \cite{GGH1}.

A different approach, based on the operations developed in \cite{GGH2} is used by M. Hern\'andez-Ketchul and J. Rodriguez-Viorato \cite{Myriam}, who wrote a Python program that is capable of computing  and printing in linear time all the distinct trivalent graphs associated to $1$-connected $2$-stratifolds up to $11$ white vertices. 

In section 2 we recall the definitions of a $2$-stratifold and its associated linear graph, providing the necessary details needed for the statement of the classification theorem for trivalent $1$-connected graphs. In section 3 we describe the general method for constructing the graphs corresponding to $1$-connected trivalent $2$-stratifolds from generating trees and skeletons, which leads to a method for counting the number of these graphs in terms of the number of black vertices of degree $3$  and the number of white vertices. In section 4 we use this approach to give explicit formulas for the case of 1 black vertex of degree $3$.  Finally  in section 5 we give a specific example to show how to compute the number of all graphs with $7$ white vertices corresponding to trivalent $1$-connected $2$-stratifolds.

\section{$2$-stratifolds and $2$-stratifold graphs.}

A  $2$-{\it stratifold} is a compact, Hausdorff space $X$ that contains a closed (possibly disconnected) $1$-manifold $X^{(1)}$ as a closed subspace with the following property: Each  point $x\in X^{(1)}$  has a neighborhood $U(x)$ homeomorphic to $CL{\times}\mathbb{R}$, where $CL$ is the open cone on $L$ for some finite set $L=\{p_1 ,\dots,p_d \}$ of cardinality $d>2$  and $X - X^{(1)}$ is a (possibly disconnected) $2$-manifold. By identifying $U(x)$ with $CL{\times}\mathbb{R}$, we call $Cp_1 {\times}\mathbb{R},\dots , Cp_d {\times}\mathbb{R}$  the {\it sheets} at $x$.\\

$X$ can be obtained as a quotient space of  a disjoint collection of circles $X^{(1)}$ and a disjoint collection $W$ of compact $2$-manifolds by attaching $W$ to $X^{(1)}$ under the attaching map $\psi$, where $\psi:\partial W\to X^{(1)}$ is a covering map such that $|\psi^{-1}(x)| >2$ for every $x\in X^{(1)}$ as in figure 1.
With suitable orientations, for a component $C$ of $\partial W$ the covering map $\psi_{|C} :C \to B\subset X^{(1)}$ is of the form $\psi (z)=z^r$, for some $r>0$.\\

We associate to a given $2$- stratifold $(X,X^{(1)} )$ an associated bi-colored labeled graph $\Gamma=\Gamma(X,X^{(1)} )$  as follows:\\

For each component $B$ of $X^{(1)}$ choose a black vertex $b$, for each component $W_i$ of $W$ choose a white vertex $w_i$, for each component $C$ of $\partial W$ choose an edge $c$. Connect $w_i$ to $b$ by the edge $c$ if $\psi (C)\subset B$.\\

We label the white vertices of the graph $\Gamma$ by assigning to $w$ the genus $g$ of $W$  (here we use Neumann's \cite{N} convention of assigning negative genus $g$ to nonorientable surfaces). We label an edge $c$ by $r$, where $r$ is the degree of the covering map $\psi_{|C} :C \to B$.\\

We say that a white vertex $w$ has genus $0$, instead of saying that the component $W$ corresponding to $w$ has genus $0$. To simplify our figures of graphs $\Gamma$, if there is no label displayed on a white vertex $w$, it is understood that the label is $0$. \\

Thus every $2$-stratifold $X$ determines uniquely a bi-colored labeled graph. Conversely, a given bi-colored labeled tree $\Gamma$ determines uniquely a $2$-stratifold $X$. \\

The association of the graph $\Gamma_X$ to the stratifold $X_\Gamma$ transforms geometrical and algebraic properties of $X_\Gamma$ into combinatorial properties of the bi-colored graph.\\

\noindent {\bf Notation}.  If $\Gamma$ is a bi-colored labeled graph corresponding to the $2$-stratifold $X$ we let $X_{\Gamma} =X$ and $\Gamma_X =\Gamma$.  An example is given in Figure 1.\\

\begin{figure}[ht]
\begin{center}
\includegraphics[width=3in]{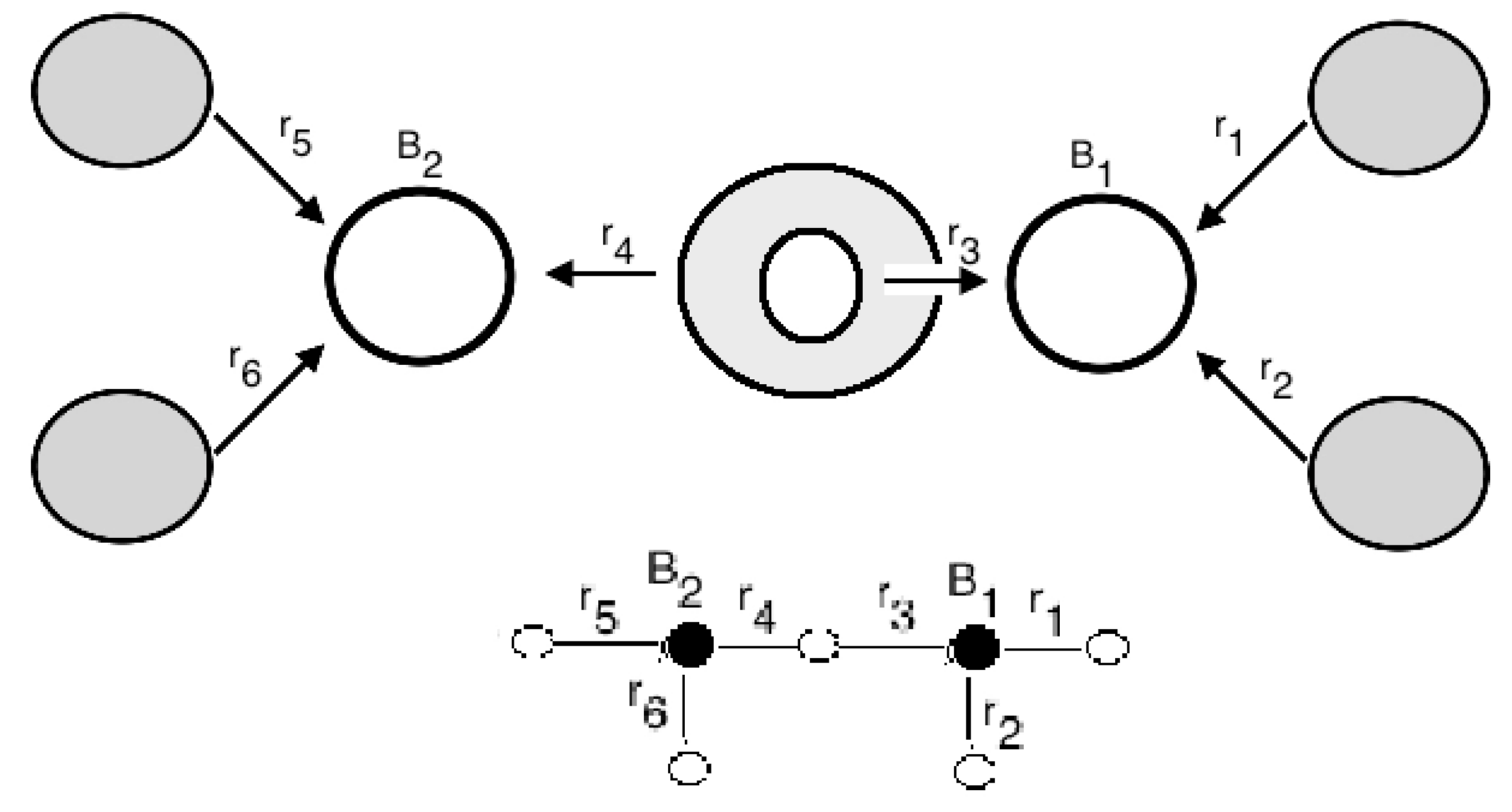} 
\end{center}
\caption{$X_\Gamma$ and $\Gamma_X$}
\end{figure}

The fundamental group $\pi_1 (X_\Gamma )$ can be computed from the bicolored graph $\Gamma_X$ (see \cite{GGH2}). In particular, if $\Gamma_X$ is a tree and all white vertices of $\Gamma_X$ have genus $0$ (i.e. correspond to punctured $2$-spheres of $X_\Gamma$), then a presentation of  $\pi_1 (X_\Gamma )$ is obtained as follows:\\
Each black vertex $b$ of $\Gamma_X$ contributes a generator, also denoted by $b$, of $\pi_1 (X)$. \\
Each white vertex $w$ incident to edges $c_1 ,\dots, c_p $ yields generators, also denoted by $c_1 ,\dots, c_p $ and a relation $c_1 \cdots c_p  =1$.\\
Each edge $c_i $ of $\Gamma_X$ between $w$ and $b$ with label $m\geq 1$ yiels a relation $b^m =c_i$.\\

The $2$-stratifold $X$ is called {\it trivalent} if every point $x\in X^{(1)}$ has a neighborhood consisting of three sheets. We do not call a $2$-manifold (i.e when $X^{(1)}=\emptyset $)  trivalent. In terms of the associated graph $\Gamma =\Gamma_X$ this means that every black vertex is incident to either one edge of label 3, or two edges one of label 1 and one of label 2, or three edges, each of label 1.

 In \cite{GGH1} we obtained a  classification theorem of simply connected trivalent $2$-stratifolds. We first review the terms used in this theorem.\\
 
(1) A {\it $(2,1)$-collapsible tree} is a bi-colored tree constructed as follows:\\
Start with a rooted tree $T$  (which may consist of only one vertex) with root $r$ (a vertex of $T$), color with white and label $0$ the vertices of $T$,  take the barycentric subdivision $sd (T)$ of $T$, color with black the new vertices (the barycenters of the edges of $T$) and finally label an edge $e$ of $sd (T)$ with $2$ (resp. $1$) if the distance from $e$ to the root $r$ is even (resp. odd). (We  allow a one-vertex tree (with white vertex) as a $(2,1)$- collapsible tree).\\

(2) The {\it reduced subgraph} $R(\Gamma)$ is defined for a bi-colored labeled tree $\Gamma$ for which the components of $\Gamma -st(\mathbb{B})$ are $(2,1)$-collapsible trees. Here $\mathbb{B}$ denotes the union of all the black vertices of degree $3$ of $\Gamma$ and $st(\mathbb{B})$ denotes the open star of $\mathbb{B}$ in $\Gamma$.  
The reduced subgraph $R(\Gamma)$ is the graph obtained from $St(\mathbb{B})$ (the closed star of $\mathbb{B}$) by attaching to each white vertex $w$ of $St(\mathbb{B})$ that is not a root, a $b12$-tree as in Figure 2, such that the terminal edge has label $2$.

\begin{figure}[ht]
\begin{center}
\includegraphics[width=3.5in]{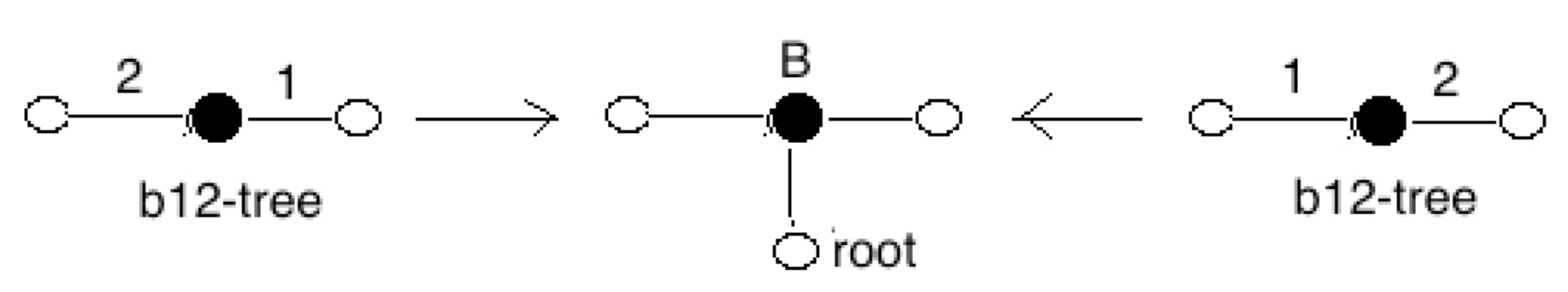} 
\end{center}
\caption{Attaching $b12$ trees}
\end{figure}

(3) A {\it horned tree} is a bi-colored tree constructed as follows:\\
Start with a tree $T$ that has at least two edges and all of whose nonterminal vertices have degree $3$.
 Color a vertex of $T$ white (resp. black) if it has degree $1$ (resp. $3$). Trisect the terminal edges of $T$ and bisect  the nonterminal edges, obtaining the graph $H_T$.
Color the additional vertices $v$ so that $H_T$ is bi-colored, that is, $v$ is colored black if $v$ is a neighbor of a terminal vertex of $H_T$ and white otherwise. Then label the edges such that every terminal edge has label $2$, every nonterminal edge has label $1$.\\

We can now state the classification theorem of \cite{GGH1}:

\begin{thm}\label{mainthm} Let $X_\Gamma$ be a trivalent $2$-stratifold with associated graph $\Gamma_X$. Let  
$\mathbb{B}$ denote the union of all the black vertices of degree $3$ of $\Gamma$ and $st(\mathbb{B})$ denote the open star of $\mathbb{B}$ in $\Gamma$.  \\
Then $X_\Gamma$ is simply connected  if and only if 
$\Gamma_X$ is a tree with all white vertices of genus $0$ and all terminal vertices white. such that the components of $\Gamma -st(\mathbb{B})$ are $(2,1)$-collapsible trees and the reduced graph $R(\Gamma)$ contains no horned tree. 
\end{thm}

\section{Skeletons}
Let $X_\Gamma$ be a $2$-stratifold whose associated graph $\Gamma_X$ has $n$ white vertices and $b$ black vertices of degree $3$. We say that $\Gamma_X$ is  trivalent $1$-connected if $X_\Gamma$ is trivalent $1$-connected.\\

We count the number of  trivalent $1$-connected graphs $\Gamma_X$ for a given number $n$ of white vertices by first counting those for a given number $b$ of black vertices of degree $3$. For such given $b$, the possible $\Gamma_X$ are obtained from the ``skeleton graphs" (defined below) that correspond to the reduced subgraphs in Theorem \ref{mainthm}.\\

\noi {\bf Generating trees}. For a given $b\geq 0$, a generating tree is an unlabeled tree with exactly $b$ black vertices and all white vertices (if any) of degree $\geq 3$. \\

\noi {\bf Skeletons}. To a generating tree $T$ we assign a skeleton $T_S$ as follows: Subdivide each edge that is incident to two black vertices and color the new vertices white. Attach edges to each black vertex such that in the resulting tree $T_S$ each black vertex has degree $3$ and all terminal vertices are white. To the white vertices $w_1 ,\dots,w_k$ of $T_S$ assign labels $T(a_1 ),\dots, T(a_k )$, where $a_i$ is an integer $\geq 1$ ($1\leq i \leq k$)\\.

Figure 3 (resp. Figure 4) shows all generating trees and their skeletons for $b=0,1,2,3$ (resp. $b=4$).\\

\begin{figure}[ht]
\begin{center}
\includegraphics[width=5in]{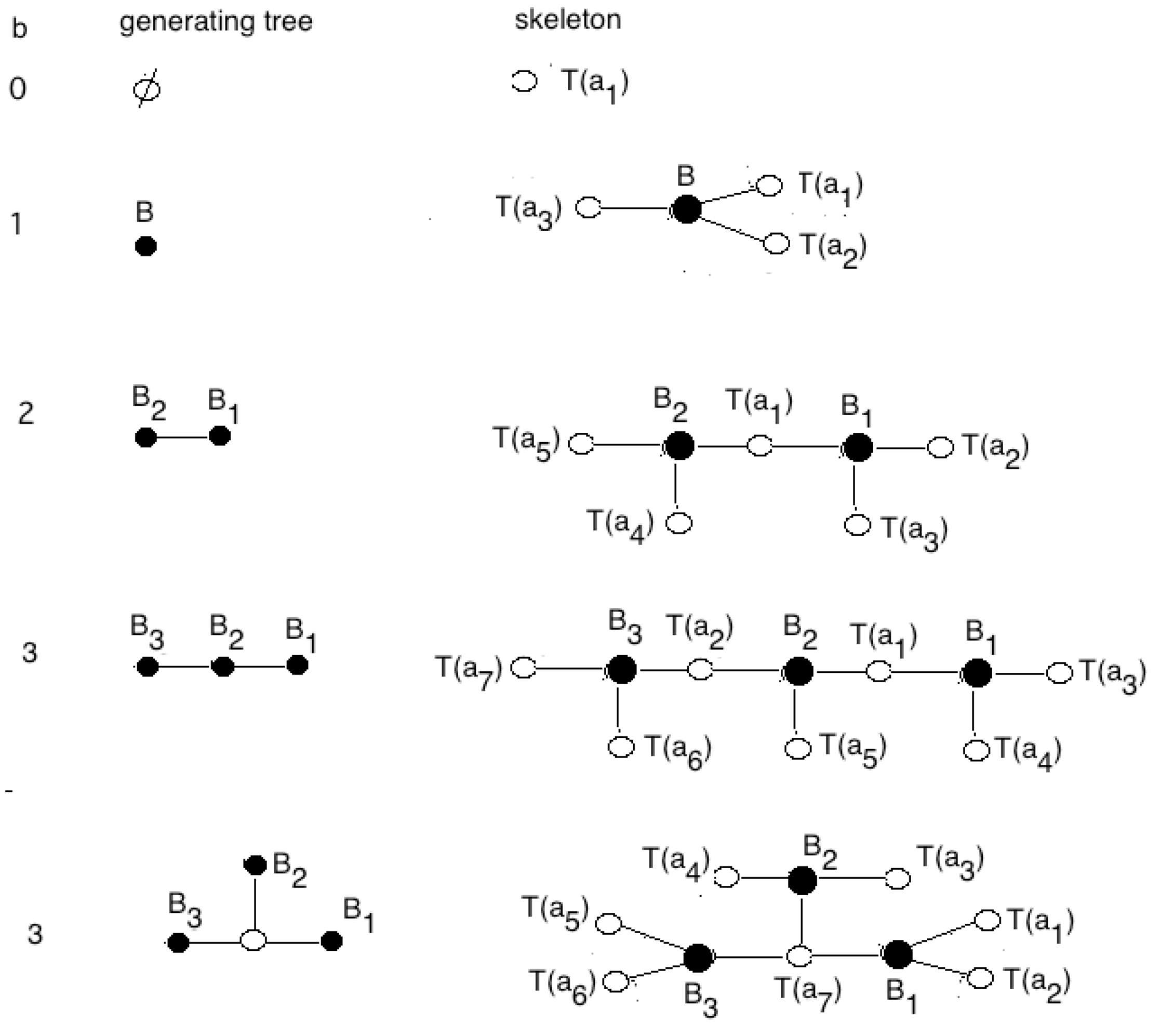} 
\end{center}
\caption{generating trees and skeletons for b=1,2,3}
\end{figure}
\begin{figure}[ht]
\begin{center}
\includegraphics[width=5in]{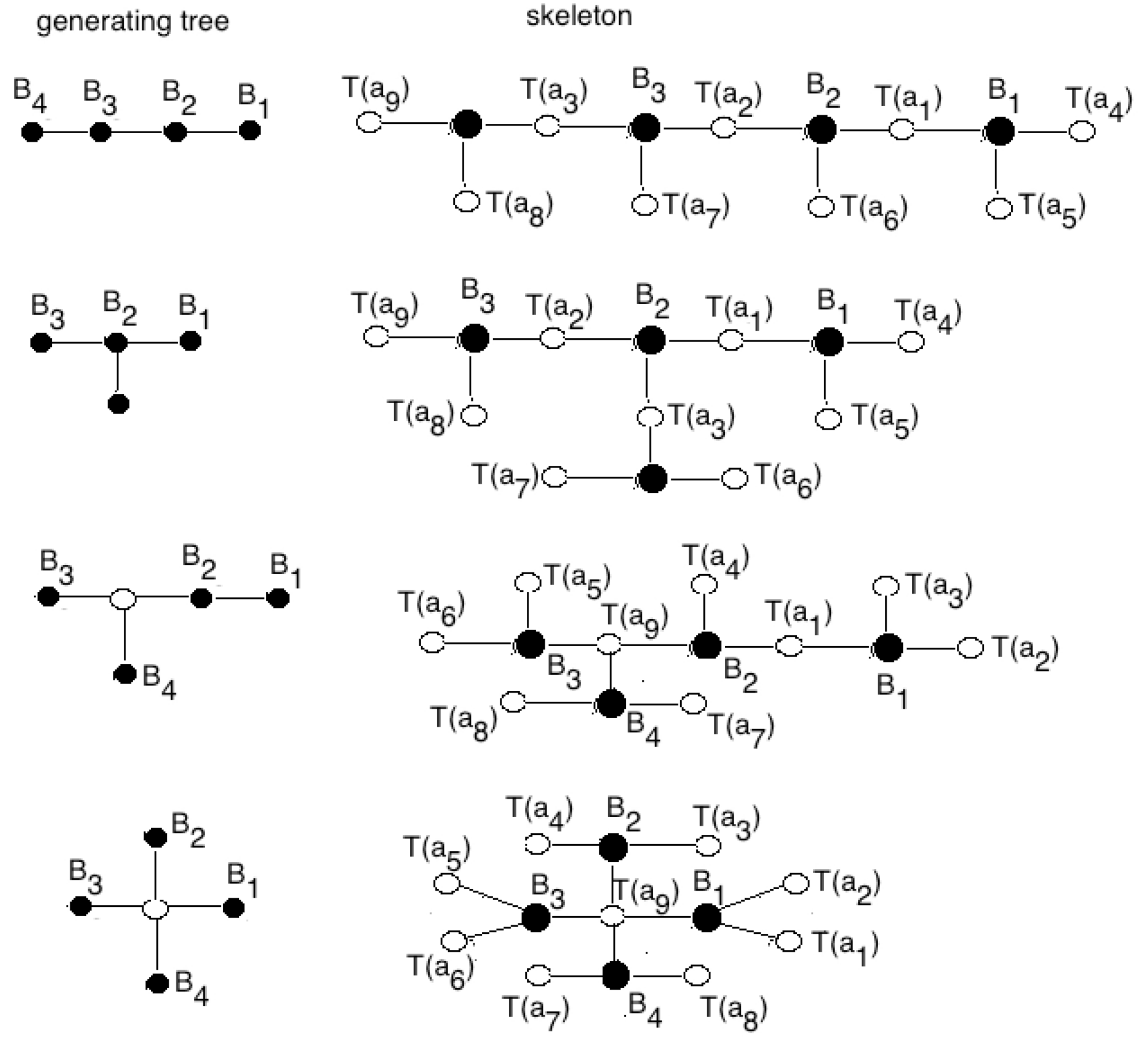} 
\end{center}
\caption{generating trees and skeletons for b=4}
\end{figure}

\noi {\bf Rooted trees}. A rooted tree $(T,r)$ is a tree $T$ with one distinguished vertex $r$, called the root of $T$.\\

\noi {\bf Bi-rooted trees}. A bi-rooted tree $(T, m, r)$ is a tree $T$ with two distinguished vertices; one called the mark $m$ and the other one called the root $r$.  We allow $m = r$, in which case one has a rooted tree.\\

\noi {\bf $d$-rooted trees}. For $d\geq 3$, a $d$-rooted tree $(T, m_1 ,\dots, m_{d-1} ; r)$ is a tree $T$ with $d$ distinguished vertices: $d-1$ marks $m_1 , \dots, m_{d-1} $ and one  root $r$. We allow $m_i =r$, for some $i$, $1\leq i \leq d-1$, but $m_i \neq m_j$ for $i\neq j$.\\

An isomorphism between bi-rooted trees $(T, m, r )$, $(T', m', r')$ (resp. $d$-rooted trees $(T, m_1 ,\dots, m_{d-1} ; r)$, $(T', m'_1 ,\dots, m'_{d-1} ; r')$) s a tree isomorphism $f:T\to T'$ such that  $f (r) = r'$ and $f (m) = m'$ (resp. $f (m_i )=m'_i $ for $i=1,\dots d-1$) \\

\begin{lem}\label{rooted} There is  a $1-1$-correspondence between $(2,1)$-collapsible trees and rooted trees. 
\end{lem}

\begin{proof} Color the vertices of the rooted tree white and bisect all edges. The new vertices are colored black. In the resulting tree $\Gamma$ assign label $2$ (resp. label $1$) to an edge that has even (resp. odd) distance to the root. Then $\Gamma$ is a $(2,1)$-collapsible tree.
\end{proof}

We now use the term {\it rooted tree} also for the associated $(2,1)$-collapsible tree.\\

By Theorem \ref{mainthm} every $1$-connected trivalent graph $\Gamma =\Gamma_X$ is obtained from $St(\mathbb{B})$ by attaching $(2,1)$-collapsible trees to the white vertices of $St(\mathbb{B})$. 

If $St(\mathbb{B})$ is connected and $\Gamma$ has $b$ black vertices of degree $3$ and $n$ white vertices, then $\Gamma$ is obtained from a skeleton (with $b$ black vertices) by attaching to each white vertex labeled $T(a_i)$ a $(2,1)$-collapsible tree having $a_i$ white vertices such that the attachment is along the mark of the corresponding bi-rooted tree. Furthermore $n=a_1 +\dots +a_k$, where $k$ is the number of white vertices of the skeleton. (If the generating tree has no white vertices, then $k=2b+1$). The symmetry group of the skeleton acts on the set of all these $\Gamma$'s and to avoid repetitions we must only count the elements in the orbits of this action. This needs to be done in such a way so that the resulting bi-colored trees do not contain horned trees. 

If $St(\mathbb{B})$ is not connected then  $\Gamma$ is obtained from a skeleton by first splitting some white non-terminal vertices. For example, if $b=2$, the skeleton splits into two cases, depending on whether $St(\mathbb{B})$ is connected or disconnected, see Figure 5. In the disconnected case the vertex of degree $2$ splits into two vertices 
and we must also consider, for a given partition $n=a_1 +a_2 +a_3 +a_4$, the number of attachments of tri-rooted trees with $a_1$ white vertices to these two vertices along two marks.

\begin{figure}[ht]
\begin{center}
\includegraphics[width=2.5in]{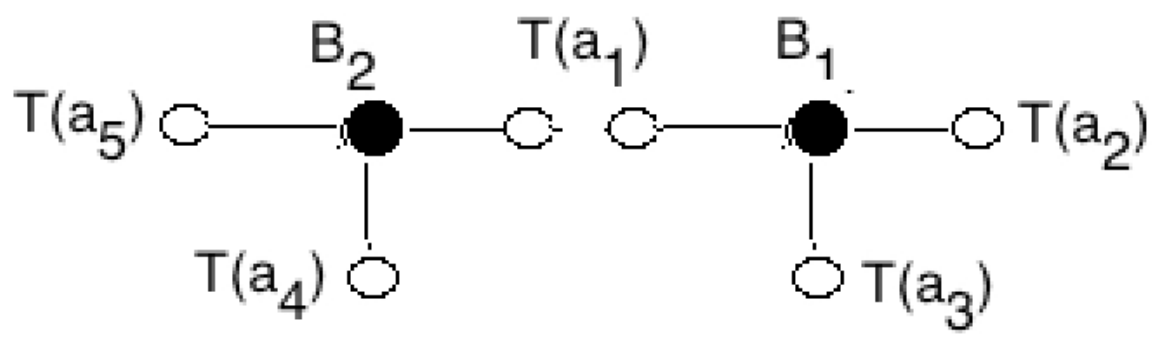} 
\end{center}
\caption{case $b=2$ disconnected}
\end{figure}

Similarly for $b\geq 2$ the skeleton splits into several cases and one must count the number of possible attachments of $k$-rooted trees for $1\leq k\leq d+1$.


\section{Number of trivalent graphs with at most one black vertex of degree $3$}

In this section we develop explicit formulas for the number of $1$-connected trivalent graphs with $n$ white vertices and one black vertex of degree $3$.

\begin{defn} $R_n$ is the number of (unlabeled) rooted trees with $n$ (white) vertices. \\
$M_a$ is the number of (isomorphism classes of) bi-rooted trees with exactly $a$ vertices.\\ 
 $U_a = M_a - R_a$  is the number of bi-rooted trees with $a$ vertices where the mark $m$ is different from the root $r$.\\
\end{defn}

The values of $R_n$ for $n\leq 30$ can be found \cite{OEIS}.\\

{\bf Case $b=0$.} Here $\Gamma_X$ is a $(2,1)$-collapsible tree. By lemma \ref{rooted} the number of distinct $1$-connected trivalent graphs $\Gamma_X$ is $R_n$.\\

{\bf Case $b=1$.} 
Here $\Gamma_X$ is obtained from a $b111$-tree (a tree with one black vertex of degree $3$ and $3$ white vertices and all edges labeled $1$) by identifying each white vertex $v_i$ of $b111$ with a white vertex of a $(2,1)$-collapsible tree $T_i$ ($i=1,2,3$) such that the reduced subgraph $R(\Gamma)$ of $\Gamma_X$ is not a horned tree. This is the case if and only if at at least one of the $v_i$'s is attached to a root of $T_i$.\\

In the skeleton graph for $b=1$ let $v_i$ be the white vertex with label $T(a_i )$. Here $T(a_i)$ is a bi-rooted tree with $a_i$ vertices and the vertex of $T(a_i )$ marked $m_i$ is identified with the vertex $v_i$ of the $b111$-graph.  The (white) edges of the bi-rooted tree $T(a_i )$ are then bisected, with the resulting vertices colored black. An edge in the bisected tree receives label $2$ (resp. $1$) if its distance to the corresponding root $r_i$ is even (resp. odd). \\

 If $\Gamma_X$ has $n$ white vertices we have $a_1 +a_2 +a_3 =n$ and in order to count all non-isomorphic  graphs with $n$ white vertices we have, by symmetry of $b111$, exactly one of the three cases $S$, $I$, $E$, below:\\

(i)  $S$ (scalene):  $a_1 > a_2 > a_3$\\
 
 (ii) $I$  (isosceles): $a_1\neq  a_2$, $a_2 = a_3$ \\
 
 (iii) $E$ (equilateral): $a_1 = a_2 = a_3$. (This occurs only when $n=3k$ for some integer $k$)\\
 
In each of the three cases let $n = a_1 +a_2 +a_3$ be a given partition. We count the number of distinct  trivalent $1$-connected graphs with $1$ black vertex of degree $3$ and $n$ white vertices.\\

(i) $S_n$: There are $M_{a_i}$  ways of attaching a birooted tree $T(a_i )$ with $a_i$  vertices to $v_i$, so there are $M_{a_1} M_{a_2} M_{a_3}$ ways of producing "scalene $(a_1 ,a_2 ,a_3 )$" trivalent trees.  However, some of these are not $1$-connected because they contain horned subtrees. So we need to subtract the number of attachments where all three vertices $v_i$ are attached to $T_i$'s along non-roots i.e. along marks $m_i$ different from the roots $r_i$. The number of these is $U_{a_1} U_{a_2} U_{a_3}$. Therefore:\\

(i) The number of distinct  trivalent $1$-connected graphs is 

$M_{a_1} M_{a_2} M_{a_3} - U_{a_1} U_{a_2} U_{a_3}$.\\

An example is shown in Figure 6 for the case $(a_1 ,a_2 ,a_3 )=(3,2,1)$.\\

\begin{figure}[ht]
\begin{center}
\includegraphics[width=5in]{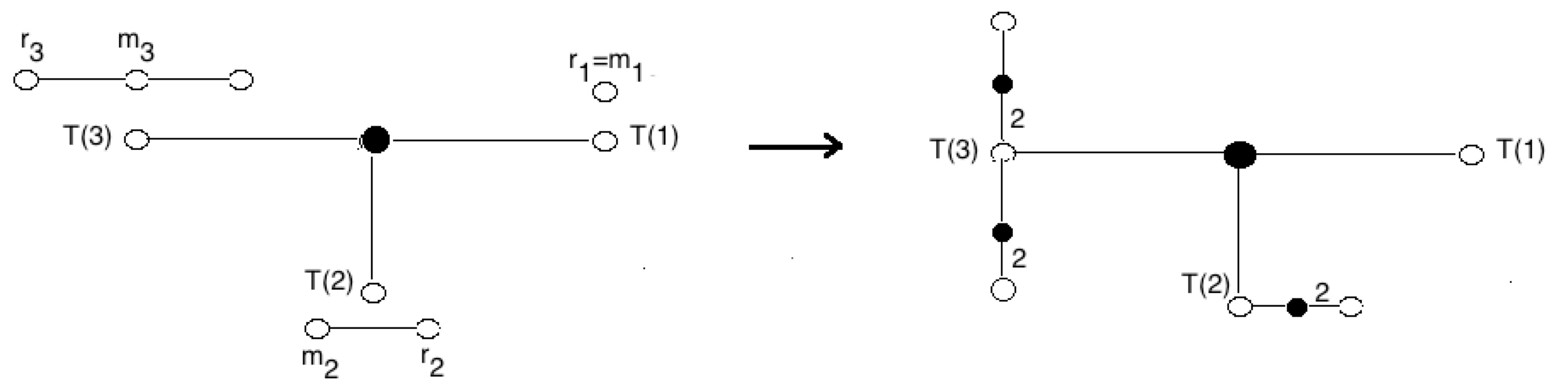}
\end{center}
\caption{Obtaining a $\Gamma_X$ from the skeleton for $b=1$}
\end{figure}

(ii) $I_n$: Let $a_1 \neq a:=a_2 =a_3$. There are $M_{a_1}$ ways to attach a birooted tree $T(a_1 )$ with $a_1$  vertices to $v_1$. Let $S_1 ,\dots, S_{M_a}$ be the distinct birooted trees with $a$ vertices. By symmetry, attaching  $S_i$ to $v_2$ and $S_j$ to $v_3$ produces the same (isomorphic) result as attaching  $S_i $ to $v_3$ and $S_j $ to $v_2$. Therefore the number of distinct graphs obtained is the number of triples $\{ T(a_1 ),S_i ,S_j \}$ with $M_a \geq i \geq j \geq 1$. To obtain the graphs corresponding to $1$-connected $2$-stratifolds we need to disregard the cases that give horned subtrees. Therefore from Lemma \ref{lem} below we obtain\\

(ii) The number of distinct isosceles trivalent $1$-connected graphs is

$M_{a_1} C(M_a +1,2)  - U_{a_1} C(U_a +1,2)$.\\

(iii) $E_n$:  Let  $a:=a_1 =a_2 =a_3$. Let $S_1 ,\dots, S_{M_a}$ be the distinct birooted trees with $a$ vertices. By symmetry, an attachment of $(S_i ,S_j ,S_k )$ to $(v_1 ,v_2 ,v_3 )$ yields isomorphic graphs if the indices $i,j,k$ are permuted. Therefore the number of distinct graphs obtained is the number of attachments of $(S_i ,S_j ,S_k )$ to $(v_1 ,v_2 ,v_3 )$  with $M_a \geq i \geq j \geq k \geq 1$. Subtracting the cases that lead to horned subtrees and using Lemma \ref{lem} we obtain:\\

(iii) The number of distinct equilateral trivalent $1$-connected graphs with $1$ black vertex of degree $3$ and $n$ white vertices is 

$\begin{cases}
   C(M_a +2,3)  - C(U_a +2,3)   & \text{if $n$ is divisible by $3$}, \\
    0  & \text{otherwise}.
\end{cases}$\\

Summing up we obtain the following Theorem.

\begin{thm}\label{thmb1} The number of distinct trivalent $1$-connected  $2$-stratifold graphs with $1$ black vertex of degree $3$ and $n$ white vertices is $S_n +I_n +E_n$.\\ 
Here $S_n=\sum   (M_{a_1} M_{a_2} M_{a_3} - U_{a_1} U_{a_2} U_{a_3})$,  where the sum is over $a_1 >a_2 >a_3$ and $a_1 +a_2 +a_3=n$\\
$I_n =\sum   (M_{a_1} C(M_a +1,2)  - U_{a_1} C(U_a +1,2))$, where the sum is over $a_1 \neq a$, $a_1 +2a =n$\\
$E_n = \begin{cases} C(M_a +2,3)  - C(U_a +2,3) & \text{ if }   3 \text{ divides } n \text{ and } 3a = n,        \\
    0   & \text{otherwise}. \end{cases}$\\
\end{thm}

\begin{lem}\label{lem} Let $m\geq 1$ and let $K=\{(k_1 ,\dots ,k_r )\in \mathbb{Z}^r \,|\,m\geq k_r \dots \geq k_2 \geq  k_1 \geq 1\}$. Then the cardinality of $K$ is $C(m+r-1,r)$.
\end{lem}

Here $C(p,q)$ is the binomial coefficient $p! /q!(p-q)! $. 

\begin{proof} An element of $K$ is a non-increasing function $k:\{1,2,\dots ,r\} \to \{1,2,\dots, m\}$, where $k(i)=k_i$. Let $\# k^{-1}(i)$ be the cardinality $k^{-1}(i)$ and denote the $m$-vector $k^{-1}=(\#k^{-1}(1), \#k^{-1}(2),\dots,\#k^{-1}(m))$ by 

\noi $\#k^{-1}(1)\,|\,\#k^{-1}(2)\,|\,\dots\,|\,\#k^{-1}(m)$ (with $m-1$ dividing bars). 

\noi From this $m$-vector delete $ \#k^{-1}(i)$ if $ \#k^{-1}(i)=0$ and replace $ \#k^{-1}(i)$ by $n$ asterisks $*$ if $ \#k^{-1}(i)=n$ to get a string of $|$'s and $*$'s.

\noi  For example if $m=8$, $r=6$ and $k=(k_1 ,\dots, k_6 )=(1,4,4,7,7,7)$, $k^{-1}=1\,|\,0\,|\,0\,|\,2\,|\,0\,|\,0\,|\,3\,|\,0 \leftrightarrow *\,|\,|\,|\,**\,|\,|\,|\,***\,|$.

This gives a bijection from the set of non-increasing functions $k:\{1,2,\dots ,r\} \to \{1,2,\dots, m\}$ to the set of all strings of length $m+r-1$ on the symbols $|$ and $*$  with exactly $r$ asterisks $*$.
\end{proof}

\section{An example for $n=7$}

In this example we show how to compute the number of $1$-connected $2$-stratifold graphs with $n=7$ white vertices. First we list a few values of $ R_n$, $M_n $, $U_n$.\\

\begin{minipage}{2in}
\begin{tabular}{c||c|c|c|}

$ n $ &  $ R_n$ &$M_n $& $U_n$\\
\hline\hline
1&1&1&0\\
2&1&2&1\\
3&2&5&3\\
4&4&13&9\\
5&9&35&26\\
6&20&95&75\\
7&48&256&208
\end{tabular}
\end{minipage}
\begin{minipage}{3in}
$R_n =$ number of rooted trees with $n$ vertices

\vspace{.05in}
$M_n =$ number of bi-rooted trees with $n$ vertices

\vspace{.05in}
$U_n =$ number of bi-rooted trees with $n$ vertices and root  different from the mark

\vspace{.05in}
\end{minipage}

\vspace{.1in} 
The table below shows how to compute the number of $1$-connected $\Gamma_X$ with $n=7$ white vertices. Here 
$b$ denotes the number of black vertices of degree $3$. The total number of non-homeomorphic $X_\Gamma$ corresponding to graphs with $n=7$ vertices is $167$.

\vspace{.1in} 
\begin{tabular}{c||c|c|c|}

&&&total number cases \\
\hline
&&&\\
$ b=0 $ & $ R_7 =48$ & & 48\\
\hline 
&&&\\
$b=1$ & $S_7 =M_4 M_2 M_1-U_4 U_2 U_1$ & & 26\\
& $I_7 =M_5 C(M_1 +1,2)-U_5 C(U_1 +1,2)$&=$35\cdot 1-25\cdot 0$&35\\
& $+M_3 C(M_2 +1,2)-U_3 C(U_2 +1,2)$&+$5\cdot 3-3\cdot 1$&12\\
& $+M_1 C(M_3 +1,2)-U_1 C(U_3 +1,2)$&+$1\cdot 15-0\cdot 6$&15\\
&$E_7 =0$&&\\
\hline
&&&\\
$b=2$& $St(\mathbb{B}$) connected: $v_0 ,v_1 ,v_2$ vertices of $St(B_1 )$&&\\
&$v_0 ,v_3 ,v_4$ vertices of $St(B_2 )$&&\\
&3 cases for middle vertex $v_0$: $a_0 =3,2,1$:&&\\
&$a_0 =3$ \quad $a_1 =a_2 =a_3 =a_4 =1$&$M_3$&5\\
&$a_0 =2$ \quad $a_1 =2$, $a_2 =a_3 =a_4 =1$&$M_2 M_2$&4\\
&$a_0 =1$ \quad $a_1 =3$, $a_2 =a_3 =a_4 =1$&$M_3 $&5\\
&\quad\qquad \quad $a_1 =2$, $a_2 =2, a_3 =a_4 =1$&$C(M_2 +1,2) $&3\\
&\quad\qquad \quad $a_1 =2$, $a_3 =2, a_2 =a_4 =1$&$C(M_2 +1,2) $&3\\
&&&\\
&$St(\mathbb{B}$) disconnected: $v_0 ,v_1 ,v_2$ vertices of $St(B_1 )$&&\\
&$v'_0 ,v_3 ,v_4$ vertices of $St(B_2 )$&& \\
&may assume tri-rooted tree is attached&& \\
&between $v_0$ and $v'_0$.  Let $a=a_0 +a'_0 \geq 2$&&\\
&&&\\
&$a=2,\,\, a_1 =2, a_3 =a_4 =1$&$M_2 M_2$&4\\
&$a=3,\,\, a_1 = a_3 =a_4 =1$&$ 5$&5\\
&&&\\
\hline
&&&\\
$b=3$& linear case &1&1\\
& star case&1&1\\
\hline
&&&\\
&Total cases for $b=0,1, 2, 3$&&167\\
\end{tabular}\\

\end{document}